\documentclass[12pt]{article}
\usepackage[dvips]{graphicx}
\usepackage{amsmath,amssymb}
\usepackage{bm}

\title{An Algorithm for the Twin Primes}
\author{}

\begin{document}
\rightline{Jan. 2014}
\rightline{~~~~}
\vskip 1cm
\centerline{\large \bf The algorithm for the $2d$ different primes}
\centerline{\large \bf and}
\centerline{\large \bf Hardy-Littlewood conjecture}
\vskip 1cm

\centerline
{{Minoru Fujimoto${}^1$
 } and
 {Kunihiko Uehara${}^2$
 }
}
\vskip 1cm
\centerline{\it ${}^1$Seika Science Research Laboratory,
Seika-cho, Kyoto 619-0237, Japan}
\centerline{\it ${}^2$Department of Physics, Tezukayama University,
Nara 631-8501, Japan}
\vskip 1cm
\centerline{\bf Abstract}
\vskip 0.2in
  We give an estimation of the existence density for the $2d$ different primes 
by using a new and simple algorithm for getting the $2d$ different primes. 
The algorithm is a kind of the sieve method, but the remainders are 
the central numbers between the $2d$ different primes. 
We may conclude that there exist infinitely many $2d$ different primes including 
the twin primes in case of $d=1$ because we can give 
the lower bounds of the existence density for the $2d$ different primes in this algorithm. 
We also discuss the Hardy-Littlewood conjecture and the Sophie Germain primes. 

\vskip 5mm
\noindent
MSC number(s): 11A41, 11Y11

\noindent
PACS number(s): 02.10.De, 07.05.Kf
\vskip 5mm

\setcounter{equation}{0}
\addtocounter{section}{0}
\section{Introduction}
\hspace{\parindent}

  Some years ago Brun\cite{Brun} dealt with the sieve methods and got 
the convergence for the inverse summation of twin primes and fixed the Brun constant, 
which corresponds to the upper bounds for the twin primes leaving their infinitude 
unsolved.\cite{McKee,Zhang} 

  The Polignac's conjecture\cite{Polignac,Shanks} involves some expectations 
about the $2d$ different prime numbers, 
which we will explain in natural way with a sieve algorithm 
where we applied the sieve of Eratosthenes. 
In addition, we discuss the Hardy-Littlewood conjecture\cite{Shanks} in the algorithm.

The basic idea of the sieve algorithm\cite{Uehara} here is that $(\text{primes\;}\pm d)$ 
are left by sifting out \\
$(\text{composite numbers\;}\pm d)$ where $d$ is any positive integer. 
As a pair of the $2d$ different prime numbers corresponds one to one to 
an arithmetic mean of these two primes, which we refer as a central number, 
the infinitude of the central number is equivalent to that of 
$2d$ different prime numbers. 

The sieve of Eratosthenes is the sieve sifting out composite numbers $np_k$ 
consequently from the natural number progression $E$ beginning at 2, 
where $n$ is integer for $n\geqq 2$ and we sift out composite numbers from $k=1$ 
in incremental order of $k$ using by the $k$-th prime number $p_k$.
When we have done the sieve using up to $p_k$, there only exist primes 
less than ${p_{k+1}}^2-1$ in the progression $E$, 
which we refer to the effective upper limit. 

\section{An algorithm for the $2d$ different primes}
\hspace{\parindent}

  The sieve for getting the central numbers of the $2d$ different primes 
is the sieve that we sift out numbers $np_k\pm d$ 
from natural number progression $N$ beginning at $2+d$, 
where $n$ is integer for $n\geqq 2$ and we sift out numbers from $k=1$ 
in incremental order of $k$ using by the $k$-th prime number $p_k$.
When we have done the sieve using up to $p_k$, 
there only exist central numbers among composite numbers 
less than ${p_{k+1}}^2-d-1$ in the progression $N$, 
which we also refer to the effective upper limit in the case. 

Especially in case of $d=1$, which corresponds to the twin primes,
\cite{Uehara,Shanks} 
only central numbers are left by the sieve. 

The reason is as follows. 
As $np_k(n\geqq 2)$ can be any composite number, we get the intersection of 
$(\text{primes}-1)$ and $(\text{primes}+1)$, namely, central numbers 
after the sieve.
So only central numbers are left in the progression $N$. 
When we look into details, we group numbers of $N$ into the primes or
the composite numbers except the central numbers or the central numbers. 
Then the primes are sifted by $2n-1$ and a composite number except the central numbers $m$
is sifted by $np_k-1$, where $np_k=m+1$ is a composite number 
because at least which of $m\pm 1$ is a composite number. 
The case that $m-1$ is a composite number is dealt as well. 
As the central numbers hold pinched by primes, 
we cannot sift them by $(\text{composite numbers}\pm 1)$ and after all
only central numbers are left in the progression $N$. 

In case of $d\ne 1$, some primes besides central numbers are left 
by the sieve, where primes are, for example, 
inner primes of prime triplets or central numbers which is primes of prime quadruplets. 
We are able to ignore these primes when we estimate the existence density 
of the $2d$ different primes because we will show later in section 5 
that the existence density of these primes cannot be a leading order, 
we have to subtract these counts though 
when we calculate the number of the $2d$ different primes.

\section{The existence density for the $2d$ different primes}
\hspace{\parindent}

  The numbers of the form $a+nb$ for any two positive coprime integers 
$a$ and $b$ form an arithmetic progression, 
where $n$ is a non-negative integer, 
and the Dirichlet prime number theorem states that this sequence 
contains infinitely many prime numbers. 
We can write the number of primes $\pi_{a,b}(x)$ up to $x$ as
\begin{eqnarray}
  \pi_{a,b}(x)
  &=&\frac{1}{\varphi(b)}Li(x)+O\left(x e^{-c_1\sqrt{\log x}}\right)\nonumber\\
  &=&\frac{1}{\varphi(b)}\left\{\frac{x}{\log x}+\frac{x}{(\log x)^2}+\frac{2x}{(\log x)^3}+O\left(\frac{x}{(\log x)^4}\right)\right\},
\end{eqnarray}
where $\varphi(b)$ is the Euler's totient function 
which gives the number of coprime to $b$ up to $b$ 
and $c_1$ is some constant. The order estimate in the second line comes from $Li(x)$. 
As this theorem is independent of $a$, 
the sieve of $np_k+a$ could be calculated by as same density as that of
Eratosthenes, which we clarify the structure of the sieve. 

When $n_1,n_2,\cdots,n_r$ are positive integers and coprime each other, 
the number of natural numbers up to $x$ which are coprime to any $n_i$ is
\begin{equation}
  S(x)=[x]-\sum_i\left[\frac{x}{n_i}\right]
          +\sum_{i<j}\left[\frac{x}{n_in_j}\right]
          -\sum_{i<j<k}\left[\frac{x}{n_in_jn_k}\right]+\cdots
          +(-1)^r\left[\frac{x}{n_1n_2\cdots n_r}\right], 
\end{equation}
where $[\ \ ]$ denotes the Gauss notation. 
This is called as the general inclusion-exclusion principle. 
This principle is applied to the number of primes up to $N$, and we get
\begin{equation}
  \pi(N)=N-\sum_i\left[\frac{N}{p_i}\right]
          +\sum_{i<j}\left[\frac{N}{p_ip_j}\right]
          -\sum_{i<j<k}\left[\frac{N}{p_ip_jp_k}\right]+\cdots
          +r-1, 
\label{e03}
\end{equation}
where $p_1,p_2,\cdots,p_r$ are primes up to $\sqrt{N}$. 
 Eq.(\ref{e03}) is correspond to counting of the primes 
by the sieve of Eratosthenes and will be the prime theorem in details 
using by the form of the zeta function. 
 When we deform it as the Euler's product 
following the Dirichlet prime number theorem, we get
\begin{equation}
  \pi(N)=cN\prod_{p\le \sqrt{N}}\left(1-\frac{1}{p}\right)+\pi(\sqrt{N})-1
        \sim cN\prod_{p\le \sqrt{N}}\left(1-\frac{1}{p}\right),
\label{e04}
\end{equation}
where $c$ is a constant determined now on and the notation $\sim$ 
in the last line means that the ratio becomes 1 for a sufficient large $N$. 
The upper limit condition for $p$ in (\ref{e04}) will be fixed 
by the Mertens' theorem\cite{Mertens} in order that we take $c=1$. 
So we put an upper limit as $p_n\le x^{e^{-\gamma}}=X(x)$ in the theorem 
\begin{equation}
  \lim_{n\to\infty}\frac{1}{\log p_n}\prod_{k=1}^n\frac{1}{1-\frac{1}{p_k}}=e^\gamma,
\end{equation}
by using 
\begin{equation}
  \prod_{p_n\le X(x)}\left(1-\frac{1}{p_n}\right)
  = \frac{1}{\log x}+O\left(\frac{1}{(\log x)^2}\right).\ \ (x\ge 4)
\end{equation}
 We write the number of central numbers up to $x$
of the $2d$ different primes as $\pi_{\text{cn}_d}(x)$,
\begin{equation}
  \pi_{\text{cn}_d}(x)\sim 2{C}_{2d}\,x\prod_{p_n\le X(x)}\left(1-\frac{2}{p_n}\right)
  = O(x) \prod_{p_n\le X(x)}\left(1-\frac{1}{p_n}\right)^2,
\label{e08}
\end{equation}
where ${C}_{2d}$ is a constant for $d$. 
The expression of the Meisel-Mertens constant $M$ 
\begin{equation}
  M=\lim_{n\to\infty}\left(\sum_{p\leqq n}\frac{1}{p}-\log(\log n)\right)
  =\gamma+\sum_p\left[\log\left(1-\frac{1}{p}\right)+\frac{1}{p}\right]
  =0.261497\cdots
\end{equation}
assures the existence for the upper limit of the coefficient $O(x)$
in the right-hand of (\ref{e08}). 

  We can estimate the Hardy-Littlewood conjecture for the twin prime
\begin{equation}
  \pi_2(x)= 2\prod_{p\ge 3}\frac{\left(1-\frac{2}{p}\right)}{\left(1-\frac{1}{p}\right)^2}
  \left\{\frac{x}{(\log x)^2}+\frac{2x}{(\log x)^3}\right\}
  +O\left(\frac{x}{(\log x)^4}\right),
\end{equation}
which is same as large as $\pi_{\text{cn}_1}(x)$, as 
\begin{equation}
  \pi_{\text{cn}_d}(x)=2C_{2d}\left\{\frac{x}{(\log x)^2}+\frac{2x}{(\log x)^3}\right\}
                      +O\left(\frac{x}{(\log x)^4}\right), 
\end{equation}
because an error estimation of the Mertens' theorem is more precise than that of 
the prime theorem. 
  We can also explain the relations in the Polignac's conjecture\cite{Polignac} 
\begin{equation}
  C_{2d}=C_2\prod_{p_k|d,\,p_k\ne 2}\frac{p_k-1}{p_k-2}
\end{equation}
in natural way, because the terms of the prime factors except 2 
in (\ref{e08}) become one series of the sieve.

\section{An algorithm for the Sophie Gelmain primes}
\hspace{\parindent}

  A prime number $p$ is a Sophie Germain prime when $2p+1$ is also prime, 
which is called a safe prime.
It has been conjectured that there are infinitely many Sophie Germain primes, 
but this remains to be proved.\cite{Shanks,Victor}

  The sieve for getting the double of the Sophie Germain primes, 
namely $(\text{the safe prime}-1)$, is the sieve sifting out numbers 
$2np_k$ and $np_k-1$ from the natural number progression $M$ beginning at $4$, 
where $n$ is integer of $n\geqq 2$ and we sift out numbers from $k=1$ 
in incremental order of $k$ using by the $k$-th prime number $p_k$.
When we have done the sieve using up to $p_k$, 
there only exist double of the Sophie Germain primes 
less than ${p_{k+1}}^2-2$ in the progression $M$, 
which we also refer to the effective upper limit in the case. 
The reason is as follows. 
As $np_k(n\geqq 2)$ can be any composite number, we get the intersection of 
$(2\times \text{primes})$ and $(\text{primes}-1)$, namely $(\text{safe primes}-1)$, 
after the sieve.
So $(2\times \text{Sophie Germain primes})$ are left in the progression $M$. 

  Here we refer an existence density for the Sophie Germain primes. 
The Dirichlet prime number theorem also holds for the sieve above, 
we get a product form of the existence density of the primes 
for the Sophie Germain primes as 
\vspace{1mm}
\begin{eqnarray}
  \frac{1}{4}\prod_{3\le p\le X(x)}\left(1-\frac{1}{p}\right)\prod_{3\le p\le X(2x)}\left(1-\frac{1}{p}\right),
\label{e12}
\end{eqnarray}
\vspace{1mm}
where the factor $\frac{1}{4}$ comes from the sieve of $p=2$ 
when only the half of even numbers are left. 
These stories are completely consistent with the Hardy-Littlewood conjecture of 
the Sophie Germain primes 
\begin{eqnarray}
  \pi_{SG}(x)&\sim&
  2\prod_{p\ge 3}\frac{\left(1-\frac{2}{p}\right)}{\left(1-\frac{1}{p}\right)^2}\frac{x}{\log x \log(2x)}\nonumber\\
  &=&\pi_2(x)\frac{\log x}{\log(2x)},
\label{e12-1}
\end{eqnarray}
\vspace{1mm}
where a sieve of $2np(p\geqq 3)$ is same as that of $np$,
because all odd numbers already sifted out by the sieve of $p=2$,
which explains why the factor 2 appears in the numerator in the product in (\ref{e12-1}).

\section{Concluding remarks}
\hspace{\parindent}

 We have discussed the $2d$ different primes or the Sophie Germain primes so far, 
but in general we can deal with a prime number $p$ when $sp+t$ is also a prime in this method, 
where $s$ and $t$ are coprimes and $2|st$. 
Moreover, as we mentioned, we deal with the primes constellation like as 
primes triplet and primes quadruplet. 

  The algorithms for the prime constellation are same as in the section 3. 
The algorithm for the central numbers of the prime triplet is the sieve sifting out $np_k-1,np_k\pm 3 \text{ or }
np_k+1,np_k\pm 3(2\le n\text{:integer})$ in order of $k$ beginning from $k=1$. 
The effective upper limit is ${p_{k+1}}^2-4$ when we sift out up to $p_k$ 
in the algorithm leaving central numbers of the prime triplet. 
In this situation, the existence density up to ${p_k}^2$ is 
\begin{equation}
  2\left(1-\frac{1}{2}\right)\left(1-\frac{2}{3}\right)\prod_{r=3}^k\left(1-\frac{3}{p_r}\right)
\label{e013}
\end{equation}
and its order will be $\displaystyle{O\left(\frac{1}{(\log k)^3}\right)}$, 
where the third power comes from the numerator in the product in (\ref{e013}). 
We note that the Hardy-Littlewood conjecture of the prime triplet 
\begin{equation}
  \pi_3(x)\sim
  9\prod_{p\ge 5}\frac{\left(1-\frac{3}{p}\right)}{\left(1-\frac{1}{p}\right)^3}\frac{x}{(\log x)^3}
  \fallingdotseq 5.716497\frac{x}{(\log x)^3} 
\label{e014}
\end{equation}
is completely reproduced here. 
The algorithm for the prime quadruplet is the sieve of $np_k\pm 2,np_k\pm 4(2\le n\text{:integer})$
in order of $k$ beginning from $k=1$. 
The effective upper limit is ${p_{k+1}}^2-5$ when we sift out up to $p_k$ 
in the algorithm leaving central numbers of the prime quadruplet. 
The existence density up to ${p_k}^2$ is 
\begin{equation}
  \left(1-\frac{1}{2}\right)\left(1-\frac{2}{3}\right)\prod_{r=3}^k\left(1-\frac{4}{p_r}\right)
\label{e015}
\end{equation}
and its order is $\displaystyle{O\left(\frac{1}{(\log k)^4}\right)}$, 
where the 4th power comes from the numerator in the product in (\ref{e015}). 
The Hardy-Littlewood conjecture of the prime quadruplet 
\begin{equation}
  \pi_4(x)\sim
  \frac{27}{2}\prod_{p\ge 5}\frac{\left(1-\frac{4}{p}\right)}{\left(1-\frac{1}{p}\right)^4}\frac{x}{(\log x)^4}
  \fallingdotseq 4.151181\frac{x}{(\log x)^4} 
\label{e016}
\end{equation}
is also reproduced again in this case. 

 It is known that the Riemann hypothesis\cite{Conrey} is equivalent to the proposition 
"there exists at least one prime between adjacent natural number squared."
As the number of twin primes up to $x$ follows to 
  $\displaystyle{\pi_2(x)=\pi_{\text{cn}_1}(x)}$,
we estimate the number of twin primes between $n^2$ and $(n+1)^2$ as
\begin{eqnarray}
  \Delta \pi_2(n^2)
  &\equiv&\pi_2((n+1)^2)-\pi_2(n^2)\nonumber\\
  &=&2C_2\left\{\frac{(n+1)^2}{(\log(n+1)^2)^2}-\frac{n^2}{(\log n^2)^2}\right\}+O\left(\frac{1}{(\log n)^3}\right)\nonumber\\
  &=&\frac{C_2}{2\log^2 n}\left\{(n+1)^2\left(\frac{\log n}{\log(n+1)}\right)-n^2\right\}+O\left(\frac{1}{(\log n)^3}\right)\nonumber\\
  &\sim& \frac{C_2}{2\log^2 n}(2n+1).
\end{eqnarray}
Therefore we could say $\Delta\pi_2(n^2)>0$ for $n>N$ with a sufficient large $N$, 
which is equivalent to the infinitude of the twin prime 
leads to the Riemann hypothesis.

\vskip 5mm

\end{document}